\documentclass{article}
\usepackage{amssymb,amsmath,theorem,euscript}

\input{epsf}

\newcounter{sec}

\def\sm{\smallskip}

\newcounter{punct}[sec]

\def\punct{\refstepcounter{punct}{\arabic{sec}.\arabic{punct}.  }}

\def\COUNTERS{\addtocounter{sec}{1}
              \setcounter{punct}{0}
          \setcounter{equation}{0}
          \setcounter{theorem}{0}
                  }

\newtheorem{theorem}{Theorem}[sec]

\begin{document}

 \def\ov{\overline}
\def\wt{\widetilde}
 \newcommand{\rk}{\mathop {\mathrm {rk}}\nolimits}
\newcommand{\Aut}{\mathop {\mathrm {Aut}}\nolimits}
\newcommand{\Out}{\mathop {\mathrm {Out}}\nolimits}
 \newcommand{\tr}{\mathop {\mathrm {tr}}\nolimits}
  \newcommand{\diag}{\mathop {\mathrm {diag}}\nolimits}
  \newcommand{\supp}{\mathop {\mathrm {supp}}\nolimits}
  \newcommand{\indef}{\mathop {\mathrm {indef}}\nolimits}
  \newcommand{\dom}{\mathop {\mathrm {dom}}\nolimits}
  \newcommand{\im}{\mathop {\mathrm {im}}\nolimits}
 
\renewcommand{\Re}{\mathop {\mathrm {Re}}\nolimits}

\def\Br{\mathrm {Br}}

\def\SL{\mathrm {SL}}
\def\SU{\mathrm {SU}}
\def\GL{\mathrm {GL}}
\def\U{\mathrm U}
\def\OO{\mathrm O}
 \def\Sp{\mathrm {Sp}}
 \def\SO{\mathrm {SO}}
\def\SOS{\mathrm {SO}^*}
 \def\Diff{\mathrm{Diff}}
 \def\Vect{\mathfrak{Vect}}
\def\PGL{\mathrm {PGL}}
\def\PU{\mathrm {PU}}
\def\PSL{\mathrm {PSL}}
\def\Symp{\mathrm{Symp}}
\def\End{\mathrm{End}}
\def\Mor{\mathrm{Mor}}
\def\Aut{\mathrm{Aut}}
 \def\PB{\mathrm{PB}}
 \def\cA{\mathcal A}
\def\cB{\mathcal B}
\def\cC{\mathcal C}
\def\cD{\mathcal D}
\def\cE{\mathcal E}
\def\cF{\mathcal F}
\def\cG{\mathcal G}
\def\cH{\mathcal H}
\def\cJ{\mathcal J}
\def\cI{\mathcal I}
\def\cK{\mathcal K}
 \def\cL{\mathcal L}
\def\cM{\mathcal M}
\def\cN{\mathcal N}
 \def\cO{\mathcal O}
\def\cP{\mathcal P}
\def\cQ{\mathcal Q}
\def\cR{\mathcal R}
\def\cS{\mathcal S}
\def\cT{\mathcal T}
\def\cU{\mathcal U}
\def\cV{\mathcal V}
 \def\cW{\mathcal W}
\def\cX{\mathcal X}
 \def\cY{\mathcal Y}
 \def\cZ{\mathcal Z}
\def\0{{\ov 0}}
 \def\1{{\ov 1}}
 \def\frA{\mathfrak A}
 \def\frB{\mathfrak B}
\def\frC{\mathfrak C}
\def\frD{\mathfrak D}
\def\frE{\mathfrak E}
\def\frF{\mathfrak F}
\def\frG{\mathfrak G}
\def\frH{\mathfrak H}
\def\frI{\mathfrak I}
 \def\frJ{\mathfrak J}
 \def\frK{\mathfrak K}
 \def\frL{\mathfrak L}
\def\frM{\mathfrak M}
 \def\frN{\mathfrak N} \def\frO{\mathfrak O} \def\frP{\mathfrak P} \def\frQ{\mathfrak Q} \def\frR{\mathfrak R}
 \def\frS{\mathfrak S} \def\frT{\mathfrak T} \def\frU{\mathfrak U} \def\frV{\mathfrak V} \def\frW{\mathfrak W}
 \def\frX{\mathfrak X} \def\frY{\mathfrak Y} \def\frZ{\mathfrak Z} \def\fra{\mathfrak a} \def\frb{\mathfrak b}
 \def\frc{\mathfrak c} \def\frd{\mathfrak d} \def\fre{\mathfrak e} \def\frf{\mathfrak f} \def\frg{\mathfrak g}
 \def\frh{\mathfrak h} \def\fri{\mathfrak i} \def\frj{\mathfrak j} \def\frk{\mathfrak k} \def\frl{\mathfrak l}
 \def\frm{\mathfrak m} \def\frn{\mathfrak n} \def\fro{\mathfrak o} \def\frp{\mathfrak p} \def\frq{\mathfrak q}
 \def\frr{\mathfrak r} \def\frs{\mathfrak s} \def\frt{\mathfrak t} \def\fru{\mathfrak u} \def\frv{\mathfrak v}
 \def\frw{\mathfrak w} \def\frx{\mathfrak x} \def\fry{\mathfrak y} \def\frz{\mathfrak z} \def\frsp{\mathfrak{sp}}
 \def\bfa{\mathbf a} \def\bfb{\mathbf b} \def\bfc{\mathbf c} \def\bfd{\mathbf d} \def\bfe{\mathbf e} \def\bff{\mathbf f}
 \def\bfg{\mathbf g} \def\bfh{\mathbf h} \def\bfi{\mathbf i} \def\bfj{\mathbf j} \def\bfk{\mathbf k} \def\bfl{\mathbf l}
 \def\bfm{\mathbf m} \def\bfn{\mathbf n} \def\bfo{\mathbf o} \def\bfp{\mathbf p} \def\bfq{\mathbf q} \def\bfr{\mathbf r}
 \def\bfs{\mathbf s} \def\bft{\mathbf t} \def\bfu{\mathbf u} \def\bfv{\mathbf v} \def\bfw{\mathbf w} \def\bfx{\mathbf x}
 \def\bfy{\mathbf y} \def\bfz{\mathbf z} \def\bfA{\mathbf A} \def\bfB{\mathbf B} \def\bfC{\mathbf C} \def\bfD{\mathbf D}
 \def\bfE{\mathbf E} \def\bfF{\mathbf F} \def\bfG{\mathbf G} \def\bfH{\mathbf H} \def\bfI{\mathbf I} \def\bfJ{\mathbf J}
 \def\bfK{\mathbf K} \def\bfL{\mathbf L} \def\bfM{\mathbf M} \def\bfN{\mathbf N} \def\bfO{\mathbf O} \def\bfP{\mathbf P}
 \def\bfQ{\mathbf Q} \def\bfR{\mathbf R} \def\bfS{\mathbf S} \def\bfT{\mathbf T} \def\bfU{\mathbf U} \def\bfV{\mathbf V}
 \def\bfW{\mathbf W} \def\bfX{\mathbf X} \def\bfY{\mathbf Y} \def\bfZ{\mathbf Z} \def\bfw{\mathbf w}
 \def\R {{\mathbb R }} \def\C {{\mathbb C }} \def\Z{{\mathbb Z}} \def\H{{\mathbb H}} \def\K{{\mathbb K}}
 \def\N{{\mathbb N}} \def\Q{{\mathbb Q}} \def\A{{\mathbb A}} \def\T{\mathbb T} \def\P{\mathbb P} \def\G{\mathbb G}
 \def\bbA{\mathbb A} \def\bbB{\mathbb B} \def\bbD{\mathbb D} \def\bbE{\mathbb E} \def\bbF{\mathbb F} \def\bbG{\mathbb G}
 \def\bbI{\mathbb I} \def\bbJ{\mathbb J} \def\bbK{\mathbb K} \def\bbL{\mathbb L} \def\bbM{\mathbb M} \def\bbN{\mathbb N} \def\bbO{\mathbb O}
 \def\bbP{\mathbb P} \def\bbQ{\mathbb Q} \def\bbS{\mathbb S} \def\bbT{\mathbb T} \def\bbU{\mathbb U} \def\bbV{\mathbb V}
 \def\bbW{\mathbb W} \def\bbX{\mathbb X} \def\bbY{\mathbb Y} \def\kappa{\varkappa} \def\epsilon{\varepsilon}
 \def\phi{\varphi} \def\le{\leqslant} \def\ge{\geqslant}

\def\UU{\bbU}
\def\Mat{\mathrm{Mat}}
\def\tto{\rightrightarrows}

\def\Gr{\mathrm{Gr}}

\def\graph{\mathrm{graph}}

\def\O{\mathrm{O}}

\def\la{\langle}
\def\ra{\rangle}

\def\B{\mathrm B}
\def\Int{\mathrm{Int}}
\def\LGr{\mathrm{LGr}}


\def\I{\mathbb I}
\def\M{\mathbb M}
\def\T{\mathbb T}
\def\S{\mathrm S}

\def\Lat{\mathrm{Lat}}
\def\LLat{\mathrm{LLat}} 
\def\Mod{\mathrm{Mod}}
\def\LMod{\mathrm{LMod}}
\def\Naz{\mathrm{Naz}}
\def\naz{\mathrm{naz}}
\def\bNaz{\mathbf{Naz}}
\def\AMod{\mathrm{AMod}}
\def\ALat{\mathrm{ALat}}

\def\Ver{\mathrm{Vert}}
\def\Bd{\mathrm{Bd}}
\def\We{\mathrm{We}}
\def\Heis{\mathrm{Heis}}

\def\bbot{{\bot\!\!\!\bot}}

\begin{center}
\Large\bf
On concentration of convolutions of double cosets at infinite-dimensional limit

\bigskip

\large\sc
Yury A. Neretin\footnote{Supported by the grants FWF, P22122, P25142.}
\end{center}

{\small For infinite-dimensional groups $G\supset K$ the double cosets $K\setminus G/K$
quite often admit a structure of a semigroup; these semigroups act
in $K$-fixed vectors of unitary representations of $G$. We show that such 
semigroups can be obtained as limits of double cosets hypergroups
(or Iwahori--Hecke type algebras) on  finite-dimensional (or finite) groups.}

\section{Introduction}

\COUNTERS

{\bf\punct An example.}
The simplest case  (Olshanski \cite{Olsh-kiado}, 1985) of the phenomenon  under discussion is the following.
Fix a positive $\alpha$.
Consider symmetric groups 
\begin{equation}
\S(\alpha)\subset \S(\alpha+1)\subset\S(\alpha+2)\subset\dots
,
\label{eq:chain}
\end{equation}
Let  assume that $\S(\beta)$ be the group of permutations of the set $\{1,2,\dots,\beta\}$.
It will be convenient to denote $G_j=\S(\alpha+j)$.
Denote by $K_j\subset G_j$ the stabilizer of points $1$, \dots, $\alpha$,
in fact $K_j$ itself is a symmetric group $\S(j)$.
We can represent elements of $G_j$ as 0-1 block  matrices of size $\alpha+(j-\alpha)$:
\begin{equation}
g=
\begin{pmatrix}
a&b\\c&d
\end{pmatrix}
\label{eq:0}
,\end{equation}
the subgroup $K_j$ consists of 0-1 matrices of the form
\begin{equation}
\begin{pmatrix}1&0\\0&d \end{pmatrix}
\label{eq:1}
.\end{equation}

Consider  two elements $g$, $h$ of $G_j$. 
 Take  large
 $N$, let us regard $g$, $h$ as elements of $G_{N}$.
 Consider the uniform probability distributions on double cosets
 $K_N g K_N$ and $K_N h K_N$. It appears, that 
  their convolution is concentrated (with probability $\to 1$ as $N\to\infty$)
  on a certain double coset 
  $$K_N r K_N \subset G_N.$$
Moreover, we can choose $r\in G_{2j}\subset G_N$ independent on $N$.

Now let $\bfG= \S(\infty)$ be the union of the chain (\ref{eq:chain}),
i.e., the group of finitely supported permutations of the countable set $\{1,2,3,\dots\}$.
Let $\bfK\subset \bfG$ be the stabilizer of $1$, \dots, $\alpha$.
Then our construction determines an associative operation $r=g\circ h$,
on $\bfK\setminus \bfG/\bfK$.

This $\circ$-product is simple and transparent.
The group $\bfG$ consists of $\alpha+\infty$ matrices of the form (\ref{eq:0}),
the subgroup $\bfK$ of matrices of the form (\ref{eq:1}).
A double coset containing $g$ is determined by the block $a$ in (\ref{eq:0}),
 and the $\circ$-product of double cosets 
 corresponds to the product of such matrices.

  It seems that mathematicians involved to infinite-dimensional classical
 groups and infinite symmetric groups believe that something in this spirit
 holds in numerous cases. However I never saw published proofs and my note covers
 this gap.
 
 \sm

{\bf\punct Hypergroups of double cosets.}
Let $G$ be a topological group, $K$ be a compact subgroup. 
Denote by $\kappa$ the probability Haar  measure on $K$.
Denote by $\delta_g$ the  delta-measure supported by
a point $g\in G$. Let $\mu*\nu$ denote the convolution
of compactly supported measures on $G$. It is determined by the condition:
for any continuous 
function $\phi$ on $G$ the following equality holds
$$
\int_G \phi(g)\,d (\mu*\nu)(g)
=\int_G\int_G \phi(g_1g_2)\,d\mu(g_1)\,d\nu(g_2)
.
$$

Denote by $K\setminus G/K$ the space of double cosets, 
recall that a {\it double coset} is a set of the form $KgK$, where $g\in G$.

For each double coset $\frg\in K\setminus G/K$ consider the measure $\sigma_\frg$
on $G$ defined by
\begin{equation}
\sigma_\frg=
\kappa*\delta_g* \kappa,\qquad \text{where $g\in\frg$}
\label{eq:sigma}
.
\end{equation}
Equivalently,  consider the map
$K\times K\to G$ given by $(k_1,k_2)\mapsto k_1gk_2$.
Then $\sigma_g$ is the pushforward of $\kappa\times\kappa$.
This measure does not depend on the choice of $g\in\frg$ and is supported by $\frg$.

A convolution of such measures can be represented as 
$$
\sigma_{\frg}*\sigma_{\frh}=
\int_{\frq\in K\setminus G/K} \sigma_\frq\,d \lambda_{\frg,\frh}(\frq)
,$$
where $\lambda_{\frg,\frh}$ is a probability measure on $K\setminus G/K$.

Thus we get an operation
$$
K\setminus G/K\,\times\, K\setminus G/K \to
\Bigl\{\text{probability measures on $K\setminus G/K$}\Bigr\}
.
$$ 
Such structures are called {\it hypergroups}
(the operation satisfies a collection of axioms, which are not important below),
see, e.g., \cite{BH}, \cite{Lit}. 

\sm

{\bf\punct Representations of hypergroups of double cosets.}
Let $\rho$ be a unitary representation of $G$ in a Hilbert space $H$.
Denote by $H^K$ the set of $K$-fixed vectors, by $P^K$ the projection
operator $H\to H^K$,
$$
P^K=\int_K \rho(k)\,d\kappa(k)
.$$
For $g\in G$ we consider  the operator
$$
P^K \rho(g) P^K=\int_{K\times K} \rho(k_1 g k_2)\,dk_1\,dk_2
$$
It depends only on the double coset $\frg$ containing $g$.
Denote by $\ov\rho(\frg):H^K\to H^K$ its restriction to $H^K$.
These operators form a representation of the hypergroup $K\setminus G/K$
in the Hilbert space $H^K$ in the following
sense:
$$
\ov\rho(\frg)\ov \rho(\frh)=\int_{K\setminus G/K}\ov\rho (\frq)\,d \lambda_{\frg,\frh}(\frq)
.
$$

{\bf \punct Dual language. Iwahori--Hecke type algebras.}
 Now let $G$ be a locally compact topological group
admitting two-side invariant Haar measure $dg$. 
 Consider the space 
$C(K\setminus G/K)$ of continuous compactly supported functions $\phi$ on 
$G$ satisfying the condition
$$
\phi(k_1gk_2)=f(g),\qquad \text{where $k_1$, $k_2\in K$}
.
$$
This space is  an algebra with respect to the convolution
$$
\phi_1*\phi_2(g)=\int_G \phi_1(h) \phi_2(h^{-1}g)\,dh
$$
Notice that $\phi(g)\,dg$ is a (sign-indefinite) measure,
and the convolution of functions corresponds to the convolution of measures.
If $\phi\in C(K\setminus G/K)$, then $\phi\,dg$ can be represented as
as
$$
\phi\,dg=\int_{\frg\in K\setminus G/K} \phi(\frg)\sigma_\frg \,dg
,$$
the measures $\sigma_\frg$ correspond to $\delta$-functions supported by 
sets $\frg$.

Several algebras $C(K\setminus G/K)$ were intensively investigated since 1950s,
namely

\sm

--- $G$ is a semisimple Lie group and $K$ is the maximal compact subgroup;
also $G$ is semisimple compact group, $K$ is a symmetric subgroup, see \cite{BG}, \cite{God}, \cite{Koo}.

\sm

--- (Hecke algebras)  $G$ is a finite Chevalley group and $K$ is the Borel subgroup, see \cite{Iwa}.

\sm

--- (affine Hecke algebras) $G$ is a $p$-adic algebraic group and $K$ is the Iwahori subgroup,
see \cite{IM}.

\sm

{\bf\punct Multiplicativity theorem.} Denote by $\U(n)$  the group
of unitary matrices of order $n$. By $\O(n)$ denote its subgroup 
consisting of real orthogonal matrices. Consider embedding
$\U(n)\to\U(n+1)$ given by
$g\mapsto \begin{pmatrix}g&0\\0&1 \end{pmatrix}$ and the direct limit
$$
\U(\infty)=\lim_{n\to\infty} \U(n)
.$$
The group $\U(\infty)$ consists of infinite unitary matrices
$g$ such that $g-1$ has only finite number of non-zero matrix elements.
In the same way we define the group $\O(\infty)\subset\U(\infty)$.

Denote by $\U(\alpha+\infty)$ the same group $\U(\infty)$ of finite unitary matrices,
but we write them as block  matrices 
$\begin{pmatrix}
a&b\\c&d
\end{pmatrix}$ of size $\alpha+\infty$.

Denote by $\O(\infty)$ the subgroup in $\U(\alpha+\infty)$ consisting of all matrices
having the following block structure
$$
h=\left(\begin{array}{ccl} 1_\alpha&0&\}\alpha\\0&u&\}\infty\end{array}\!\!\!\!\!\!\!\!\!\!\!\!\!\!\!\!\right)
\qquad, \quad\text{where $u\in\O(\infty)$}
$$
and $1_\alpha$ is the unit matrix of size $\alpha$.
 Now set
$$
\bfG=\U(\alpha+\infty),\quad \bfK=\O(\infty).
$$

Consider double cosets 
$\bfK\setminus \bfG/\bfK$, i.e. finite unitary matrices
defined up to the equivalence
$$
\begin{pmatrix}
a&b\\c&d
\end{pmatrix}
\sim
\begin{pmatrix}
1_\alpha&0\\0&u
\end{pmatrix}
\begin{pmatrix}
a&b\\c&d
\end{pmatrix}
\begin{pmatrix}
1_\alpha&0\\0&v
\end{pmatrix},\quad\text{where $u$, $v\in\O(\infty)$.}
$$
There is no Haar measure on $\bfK$ and therefore we can not define canonical
measures on $\bfK g\bfK$.

Consider a unitary representation $\rho$ of the group $\bfG$ in a Hilbert space
$H$. As above consider the space $H^\bfK$ of $\bfK$-fixed vectors and  the projection operator
$H\to H^\bfK$. Again, we define an operator 
$P^\bfK\rho(g)P^\bfK$ and the operator 
$$
\ov\rho(\frg):=P^\bfK\rho(g)P^\bfK\Bigr|_{H^\bfK}
.$$

\begin{theorem}
\label{th:multi-1}
{\rm (Olshanski, \cite{Olsh-GB})}
{\rm a)} For any $\frg$, $\frh\in \bfK\setminus \bfG/\bfK$
there exists an element $\frg\circ\frh\in \bfK\setminus \bfG/\bfK$ 
such that for any unitary representation $\rho$ of $\bfG$,
$$
\ov\rho(\frg)\,\ov\rho(\frh)=\ov\rho(\frg\circ \frh)
.
$$

{\rm b)} Moreover, the operation $\frg\circ\frh$
is given by the formula 
\begin{multline}
\begin{pmatrix}
a&b\\c&d
\end{pmatrix}
\circ
\begin{pmatrix}
p&q\\r&t
\end{pmatrix}:=
\begin{pmatrix}
a&b&0\\c&d&0\\0&0&1_\infty
\end{pmatrix}
\begin{pmatrix}
1_\alpha&0&0\\
0&0&1_\infty\\
0&1_\infty&0
\end{pmatrix}
\begin{pmatrix}
p&q&0\\r&t&0\\0&0&1_\infty
\end{pmatrix}
=\\=
\begin{pmatrix}
ap&aq&b\\
cp&cq&d\\
r&t&0
\end{pmatrix}
\sim
\begin{pmatrix}
ap&b&aq\\
cp&d&cq\\
r&0&t
\end{pmatrix}
.
\end{multline}
\end{theorem}

A detailed description of this semigroup is contained
in \cite{Ner-book},  Section IX.4.

\sm

{\bf\punct The purpose of the paper.} We wish to show that in a 
certain sense the semigroup
$\O(\infty)\setminus \U(\alpha+\infty)/\O(\infty)$
is the limit of hypergroups $\O(n)\setminus \U(\alpha+n)/\O(n)$
as $n\to\infty$. The precise statement is in the next section.

First examples of semigroups of double cosets were discovered
in 70s by Ismagilov and Olshanski, see \cite{Ism}, \cite{Olsh-kiado},
\cite{Olsh-tree}.
Now lot of semigroups $\bfK\setminus\bfG/\bfK$ are known,
usually they admit explicit realizations, see  \cite{Olsh-GB},
\cite{Olsh-symm}, \cite{Ner-book}--\cite{Ner-invariant}.

Below we prove several well-representative 'limit theorems'
for semigroups $\bfK\setminus\bfG/\bfK$ for real classical
and symmetric groups. It seems that for $p$-adic groups
\cite{Ner-p}
this approach fails.


\section{Double cosets. Classical groups}

\COUNTERS

{\bf \punct The limit theorem for $\bfG=\U(\alpha+\infty)$, $\bfK=\O(\infty)$.} Consider a group $G_N:=\U(\alpha+k+N)$. We equip $G_N$ with the metric
$\|g-h\|$, where $\|\cdot\|$ denotes the norm of an operator in a Euclidean
space.
Consider the subgroup $\U(\alpha+k)\subset G_N$ consisting of matrices
$$
\left(
\begin{matrix}
a&b&0&\} \alpha\\c&d&0&\}k\\0&0&1_N&\}N
\end{matrix}
\!\!\!\!\!\!\!\!\!\!\!\!\!\!\right)
\qquad.$$

Consider the subgroup $K_N:=\O(k+N)\subset G_N$ consisting of real matrices
$$
\begin{pmatrix}
1_\alpha&0&0\\
0&u&v\\
0&w&y
\end{pmatrix}
.$$
Denote by $\kappa_N$ the probability Haar measure supported by $K_N$.

Let $g$, $h\in\U(\alpha+k)$. We define their $\circ_N$-product,
$$
g\circ_N h\in K_N\setminus G_N/K_N
$$
 by
\begin{multline*}
\begin{pmatrix}
a&b&0\\c&d&0\\0&0&1_N
\end{pmatrix}
\circ_N
\begin{pmatrix}
p&q&0\\r&t&0\\0&0&1_N
\end{pmatrix}
:=\\:=
\begin{pmatrix}
a&b&0&0\\c&d&0&0\\0&0&1_k&0
\\0&0&0&1_{N-k}
\end{pmatrix}
\begin{pmatrix}
1_\alpha&0&0&0\\
0&0&1_k&0\\
0&1_k&0&0
\\0&0&0&1_{N-k}
\end{pmatrix}
\begin{pmatrix}
p&q&0&0\\r&t&0&0\\0&0&1_k&0\\0&0&0&1_{N-k}
\end{pmatrix}
=\\=
\begin{pmatrix}
ap&aq&b&0\\
cp&cq&d&0\\
r&t&0&0
\\0&0&0&1_{N-k}
\end{pmatrix}
\end{multline*}
and we take the double coset containing the latter matrix.

On the other hand, for any two elements $g$, $h\in\U(\alpha+k)$  we consider the measure
$\tau_{g,h}$ on $\U(\alpha+k+N)$ given by 
\begin{equation}
\tau_{g,h}:=\kappa_N * \delta_g* \kappa_N *\delta_h*\kappa_N
.
\label{eq:tau-1}
\end{equation}
In notation (\ref{eq:sigma}),
\begin{equation}
\tau_{g,h}=\sigma_g*\sigma_h
\label{eq:tau-2}
\end{equation}

\begin{theorem}
\label{th:1}
Fix $\alpha$. For a given $k$ and any $\epsilon>0$, $\delta>0$
there exists $N_0$ such that for $N\ge N_0$ for any $g$, $h\in\U(\alpha+k)$
the measure $\tau_{g,h}$ of the $\epsilon$-neighborhood of the coset
$g\circ_N h$
is $> 1-\delta$.
 \end{theorem}
 
 {\sc Proof.}
 Consider the following  measure on $G_N$:
 $$
\wt\tau_{g,h}:=\delta_g* \kappa_N *\delta_h
.$$
Consider a subset $A\subset \U(\alpha+k+N)$ invariant with respect to left and right
translations by elements of $\O(k+N)$. Then 
$$
\wt\tau_{g,h}(A)=\tau_{g,h}(A)
.
$$
 For this reason we will estimate  the measure $\wt\tau_{g,h}$
 of the $\epsilon$-neighborhood of $g\circ_N h$,
  \begin{equation}
\begin{pmatrix}
a&b&0\\c&d&0\\0&0&1_N
\end{pmatrix}
\begin{pmatrix}
1_\alpha&0&0\\
0&u&v\\
0&w&y
\end{pmatrix}
\begin{pmatrix}
p&q&0\\r&t&0\\0&0&1_N
\end{pmatrix}
=\qquad\qquad\qquad\qquad\qquad\qquad
\label{eq:triple-0}
\end{equation}
\begin{equation}
\qquad\qquad\qquad\qquad\qquad\qquad
=
\begin{pmatrix}
ap+\boxed{bur}&aq+\boxed{but}&bv\\
cp+\boxed{dur}&cq+\boxed{dut}&dv\\
wr&wt&y
\end{pmatrix}
\label{eq:triple-1}
\end{equation}
The matrix  $\begin{pmatrix}u&v\\w&y \end{pmatrix}$ is orthogonal.
Therefore, for each $i$
$$
\sum_{j=1}^k u_{ij}^2+\sum_{l=1}^N v_{il}^2=1
$$
If $N$ is large, then a given summand of this sum is $\simeq 0$ with probability near 1.
The size of the matrix $u$ is fixed, the size of $\begin{pmatrix}u&v\\w&y \end{pmatrix}$
tends to infinity.
Therefore, for any $\delta>0$ and any $\epsilon'>0$ we can choose $N$ such that
$\|u\|<\epsilon'$ on the set of measure $>1-\delta$. This is our main argument.
We wish to show that if $\|u\|$ is small, then (\ref{eq:triple-1})
is contained in a small neighborhood of the set $g\circ_N h$.

First, boxed terms in (\ref{eq:triple-1})  are small.

Next, we can multiply the middle factor in (\ref{eq:triple-0})
 by elements of $\O(N)$ from left and right.
  Such transformations do not change a double coset containing
(\ref{eq:triple-1}),
\begin{equation}
\begin{pmatrix}
1_\alpha&0&0\\
0&u&v\\
0&w&y
\end{pmatrix} \mapsto
\begin{pmatrix}
1_\alpha&0&0\\
0&u&v\eta\\
0&\xi w&\xi y\eta
\end{pmatrix}, \qquad \xi,\eta\in\O(N)
\label{eq:middle-1}
.
\end{equation}

We have $u^t u+w^tw=1_k$. Since $\|u\|$ is small%
\footnote{In particular, $N>k$, otherwise $\|u\|=1$.}%
, the matrix $w$ is 'almost isometry'
$\R^k\to \R^{N}$. Therefore it can be reduced by multiplication
$w\mapsto \xi w$ to the form
$$w=\begin{pmatrix}1_k+\tau\\0\end{pmatrix},\quad  \text{where $\|\tau\|$ is small}.$$
In the same way, we can reduce $v$ to that form
$$
v=\begin{pmatrix}1_k+\sigma&0\end{pmatrix}
,\quad  \text{where $\|\sigma\|$ is small}.
$$ 
Thus (\ref{eq:middle-1}) is reduced to the form
$$
\begin{pmatrix}
1_\alpha&0&0&0\\
0&u&1_k+\sigma&0\\
0&1_k+\tau&y_{11}&y_{12}\\
0&0&y_{21}&y_{22}
\end{pmatrix}
.
$$
This matrix is orthogonal, therefore $\|y_{11}\|$, $\|y_{12}\|$, $\|y_{21}\|$
are small. Referring to orthogonality again, we get that $y_{22}=z+\phi$,
where $z\in\O(N-k)$ and $\|\phi\|$ is small. 
We also can transform $y_{22}\to \gamma y_{22} \delta$, where $\gamma$, $\delta\in \O(N-k)$.
Such transformations do not change a double coset containing
(\ref{eq:triple-1}). Therefore, we can assume $y_{22}=1+\psi$, where $\|\psi\|$ is small.
Thus we can replace
$$
\begin{pmatrix}
1_\alpha&0&0\\
0&u&v\\
0&w&y
\end{pmatrix}
\to \begin{pmatrix}
1_\alpha&0&0&0\\
0&0&1_k&0\\
0&1_k&0&0\\
0&0&0&1_{N-k}
\end{pmatrix}\,+ \Bigl\{\text{matrix with small norm}\Bigr\}
$$
without changing of  double coset (\ref{eq:triple-1}).
But after this replacement  we get a matrix closed to $g\circ_N h$.
\hfill $\square$

\sm

{\bf\punct A more complicated case.}
Now let $\bfG=\U(\alpha+m\infty)$. It is the same group $\U(\infty)$,
but we write matrices in the block form
\begin{equation}
g=\begin{pmatrix}
a&b_1&\dots &b_m\\
c_1&d_{11}&\dots&d_{1m}\\
\vdots&\vdots&\ddots&\vdots\\
c_m&d_{m1}&\dots&d_{mm}
\end{pmatrix}
\label{eq:g-big}
\end{equation}
of size $\alpha+\infty+\dots+\infty$.
Consider its subgroup $\bfK\simeq\O(\infty)$
consisting of matrices 
\begin{equation}
\begin{pmatrix}
1_\alpha&0&0&\dots&0\\
0&u&0&\dots&0\\
0&0&u&\dots&0\\
\vdots&\vdots&\vdots&\ddots&\vdots\\
0&0&0&\dots&u
\end{pmatrix}
,
\label{eq:u-big}
\end{equation}
where $u\in \O(\infty)$.

 {\it Double cosets $\bfK\setminus \bfG/\bfK$ admit a natural structure of a semigroup
and the multiplicativity theorem for the pair $\bfG\supset \bfK$ also holds}
(\cite{Ner-faa}). To simplify notation, we set $m=2$(a pass from 2 to 
arbitrary $m$ is straightforward).

First, we define the $\circ$-multiplication,
\begin{multline}
\begin{pmatrix}
a&b_1&b_2\\
c_1&d_{11}&d_{12}\\
c_2&d_{21}&d_{22}
\end{pmatrix}
\circ
\begin{pmatrix}
p&q_1&q_2\\
r_1&t_{11}&t_{12}\\
r_2&t_{21}&t_{22}
\end{pmatrix}
=\\=
\scriptsize
\begin{pmatrix}
a&b_1&0&b_2&0\\
c_1&d_{11}&0&d_{12}&0\\
0&0&1_\infty&0&0\\
c_2&d_{21}&0&d_{22}&0\\
0&0&0&0&1_\infty
\end{pmatrix}
\begin{pmatrix}
1_\alpha&0&0&0&0\\
0&0&1_\infty&0&0\\
0&1_\infty&0&0&0\\
0&0&0&0&1_\infty\\
0&0&0&1_\infty&0
\end{pmatrix}
\begin{pmatrix}
p&q_1&0&q_2&0\\
r_1&t_{11}&0&t_{12}&0\\
0&0&1_\infty&0&0\\
r_2&t_{21}&0&t_{22}&0\\
0&0&0&0&1_\infty
\end{pmatrix}
=\\=
\normalsize
\begin{pmatrix}
ap& |&  aq_1 & b_1&&  aq_2 & b_2
\\
-& + & -& -& -&-
\\
c_1 p&|&  c_1 q_1& d_{11} && c_1 q_2& d_{12} 
\\
r_1 &|&  t_{11}&0 &&  t_{12}&0 
\\
&|&&&&
\\
c_2 p&|&  c_2 q_1& d_{21} && c_2 q_2& d_{22} 
\\
r_2 &|&  t_{21}&0 &&  t_{22}&0 
\end{pmatrix}
\label{eq:product-big}
\end{multline}

Next, we introduce an operation $\circ_N$ as follows.
Consider the group $G_N:=\U(\alpha+ 2(k+N))$. We write its elements as
block matrices of size $\alpha+k+N+k+N$. Sometimes we will subdivide 
these matrices and write them as block matrices of size
$\alpha+k+k+(N-k)+k+k+(N-k)$. Consider the subgroup $\U(\alpha+2k)\subset G_N$
consisting of matrices
\begin{equation}
g=
\begin{pmatrix}
a&b_1&0&b_2&0\\
c_1&d_{11}&0&d_{12}&0\\
0&0&1_N&0&0\\
c_2&d_{21}&0&d_{22}&0\\
0&0&0&0&1_N
\end{pmatrix}
\label{eq:g}
\end{equation}
and the subgroup $K_N\simeq\O(k+N)\subset G_N$ consisting of matrices
\begin{equation}
\begin{pmatrix}
1_\alpha&0&0&0&0\\
0&u&v&0&0\\
0&w&t&0&0\\
0&0&0&u&v\\
0&0&0&w&t
\end{pmatrix}, \qquad \text{where $\begin{pmatrix}u&v\\w&t\end{pmatrix}\in \O(k+N)$.}
\label{eq:double-O}
\end{equation}
Denote by $\kappa_N$  the  Haar measure of $K_N$ regarded as a measure on $G_N$.
For $g$, $h\in\U(\alpha+2k)$ denote by $\tau_{g,h}$
the measure on $G_N$ given by
$$
\tau_{h,h}:=\kappa_N*\delta_g*\kappa_N*\delta_h*\kappa_N
$$
as above.

Denote by $J$ the following matrix
\begin{equation}
J_N:=
\begin{pmatrix}
1_\alpha&0&0&0&0&0&0\\
0&0&1_k&0&0&0&0\\
0&1_k&0&0&0&0&0\\
0&0&0&1_{N-k}&0&0&0\\
0&0&0&0&0&1_k&0\\
0&0&0&0&1_k&0&0\\
0&0&0&0&0&0&1_{N-k}
\end{pmatrix}\in \U(\alpha+2(k+N))
\label{eq:JN-big}
\end{equation}
For $g$, $h\in\U(\alpha+2k)$ we define
$$g\circ_N h\in K_N\setminus G_N/K_N$$
 by
\begin{equation}
g\circ_N h=K_N \cdot gJ_Nh \cdot K_N,\qquad\text{where $J_N$ is given by (\ref{eq:JN-big})}
.
\label{eq:gJh}
\end{equation}

\begin{theorem}
\label{th:2}
Fix $\alpha$. For any $k$ for any $\epsilon>0$, $\delta>0$
there exists $N_0$ such that for any $N\ge N_0$ the measure $\tau_{g,h}$
of the $\epsilon$-neighborhood of the coset
$g\circ_N h$ 
is $\ge 1-\delta$.
\end{theorem}

{\sc Proof} repeats the proof of Theorem \ref{th:1}.
We evaluate analog of the product (\ref{eq:triple-1}), i.e.
a product 
\begin{equation}
gUh
,\end{equation}
where $g$ is given by (\ref{eq:g}), $U$ is given by (\ref{eq:double-O}),
and 
\begin{equation}
h=\begin{pmatrix}
p&q_1&0&q_2&0\\
r_1&t_{11}&0&t_{12}&0\\
0&0&1_N&0&0\\
r_2&t_{21}&0&t_{22}&0\\
0&0&0&0&1_N
\end{pmatrix}
\label{eq:triple-2}
.
\end{equation}
Again we note that $\|u\|$ is small on a subset whose complement has small measure.
Repeating one-to-one the same steps,
we change $U$ in (\ref{eq:triple-2}) by
$$
J_N+ \Bigl\{\text{matrix with small norm}\Bigr\}
$$
and come to the desired statement.
\hfill $\square$


\section{Hypergroup of conjugacy classes and  operator colligations}

\COUNTERS

{\bf\punct Hypergroups of conjugacy classes.} Consider a topological  group
$G$ and a compact subgroup $K$. Denote by $G//K$ the conjugacy classes
on $G$ with respect to $K$, i.e., the quotient of $G$ with respect to the equivalence relation
$$
g\sim hgh^{-1}, \qquad \text{where $h$ ranges in $K$.}
$$
 Consider the probability Haar measure
$\kappa$ on $K$. For each $g\in G$ consider the map
$K\to G$ given by $h\mapsto hgh^{-1}$. Denote by $\nu_g$ the image of 
$\kappa$ under this map. It is readily seen that $\nu_g$
depends only on the $K$-conjugacy class $\frg$ containing $g$. Again, 
measures $\nu_\frg$ form a hypergroup, i.e.,
$$
\nu_\frg*\nu_\frq=\int_{G//K} \nu_\frr\,d\psi_{\frg,\frq}(\frr),
$$
where $d\psi_{\frg,\frq}$ is a probability measure on $G//K$

\sm

{\bf\punct Operator colligations.}
Now consider the group $\bfG=\U(\alpha+\infty)$
and its subgroup $\bfK=\U(\infty)$.
{\it Operator colligations}%
\footnote{A synonym: node.} are conjugacy classes 
 $\bfG//\bfK$.
There is a well-defined associative operation
$$
\bfG//\bfK\times \bfG//\bfK\,\to\, \bfG//\bfK
$$
given by
\begin{equation}
\begin{pmatrix}
a&b\\c&d
\end{pmatrix}
\circ
\begin{pmatrix}
p&q\\r&t
\end{pmatrix}:=
\begin{pmatrix}
a&b&0\\c&d&0\\0&0&1_\infty
\end{pmatrix}
\begin{pmatrix}
p&0&q\\0&1_\infty&0 \\r&0&t
\end{pmatrix}
=
\begin{pmatrix}
ap&b&aq\\
cp&d&cq\\
r&0&t
\end{pmatrix}
.
\end{equation}

This semigroup (with several modifications) is a classical topic of operator theory
and system theory, see, e.g.,
\cite{Bro}, \cite{Dym}. 

We wish to show that the semigroup
$\U(\alpha+\infty)//\U(\infty)$ is a limit  of hypergroups
$\U(\alpha+n)//\U(n)$ as $n\to\infty$.

\sm

{\bf\punct The limit theorem.} Again, consider the group $G_N=\U(\alpha+k+N)$ and its
subgroups $\U(\alpha+k)$ and $K_N=\U(k+N)$. For $g$, $h\in \U(\alpha+k)$ 
consider the measure
$$
\tau_{g,h}:=\nu_\frg*\nu_\frh
$$
on $G_N$. 

On the other, for $g$, $h\in \U(\alpha+k)$ we define
the element
$$
g J h J^{-1}
,$$
where
$$
J:=\begin{pmatrix}
1_\alpha&0&0&0\\
0&0&1_k&0\\
0&1_k&0&0\\
0&0&0&1_{N-k}
\end{pmatrix}
$$
 and denote by
 $g\circ_N h$
the corresponding conjugacy class $\in\bfG//\bfK$.

\begin{theorem}
Fix $\alpha$. For a given $k$ for any $\epsilon>0$, $\delta>0$
there exists $N_0$ such that for all $N\ge N_0$ the measure
$\tau_{g,h}$ of the $\epsilon$-neighborhood 
of $g\circ_N h$ is $>1-\delta$.
\end{theorem}

{\sc Proof.}
Denote by $\wt\tau_{g,h}$ the image of the measure $\kappa_N$ under the map 
$$
x\mapsto g x h x^{-1},\qquad \text{where $x$ ranges in $K_N$}
.$$
The measure $\tau_{g,h}$ coincides with the image of the measure
$\kappa_N\times \kappa_N$
under the map
$$
(x,z)\mapsto  z g x h x^{-1}z^{-1},\qquad \text{where $x$, $z$ range in $K_N$}
$$
For any set $L\subset G_N$ invariant with respect to $K_N$-conjugations, we have
$\wt\tau_{g,h}(L)=\tau_{g,h}(L)$. Next, we evaluate
\begin{equation}
gxhx^{-1}
\label{eq:gxh}
\end{equation}
 for $x\in K_N$:
\begin{multline*}
\begin{pmatrix}
a&b&0\\c&d&0\\0&0&1_N
\end{pmatrix}
\begin{pmatrix}
1_\alpha&0&0\\
0&u&v\\
0&w&y
\end{pmatrix}
\begin{pmatrix}
p&q&0\\r&t&0\\0&0&1_N
\end{pmatrix}
\begin{pmatrix}
1_\alpha&0&0\\
0&u^*&w^*\\
0&v^*&y^*
\end{pmatrix}
=\\=
\begin{pmatrix}
ap+bur&aqu^*+butu^*+bvv^*&aqw^*+butw^*+bvy^*\\
cp+dur&cqu^*+dutu^*+dvv^*&cqw^*+dutw^*+dvy^*\\
wr&wtu^*+yv^*&wtw^*+yy^*
\end{pmatrix}
\end{multline*}
Again, we can assume that $\|u\|$ is small, therefore the last matrix can be represented as
\begin{equation}
\begin{pmatrix}
ap &bvv^*& aqw^*+bvy^*\\
cp& dvv^*& cqw^*+dvy^*\\
wr&  yv^*& wtw^*+yy^*
\end{pmatrix}+\Bigl\{\text{small matrix}\Bigr\}
\label{eq:prod-conj-2}
.
\end{equation}
Next, we have
$$
\begin{pmatrix}1_k&0\\0&1_N\end{pmatrix}
=
\begin{pmatrix}u&v\\w&t\end{pmatrix}
\begin{pmatrix}u^*&w^*\\v^*&t^*\end{pmatrix}=
\begin{pmatrix}
uu^*+vv^*& uw^*+vy^*\\
wu^*+y v^*&ww^*+vv^*
\end{pmatrix}
.
$$
Since $\|u\|$ is small, we have 
$$
vy^*\simeq 0,\quad yv^*\simeq 0,\quad vv^*\simeq 1_k.
$$
Therefore (\ref{eq:prod-conj-2}) can be written as 
\begin{equation}
\begin{pmatrix}
ap &b& aqw^*\\
cp& d& cqw^*\\
wr& 0& wtw^*+yy^*
\end{pmatrix}+\Bigl\{\text{small matrix}\Bigr\}
\label{eq:prod-conj-3}
.
\end{equation}
Next, we can conjugate a matrix $x\in\U(k+N)$ in
(\ref{eq:gxh})
by elements of $\U(N)$, this leads to a conjugation of the whole expression
(\ref{eq:gxh}). In this way, we can reduce $w$ to the form
$$
w=\begin{pmatrix} 1_k\\0\end{pmatrix}+\Bigl\{\text{small matrix}\Bigr\}
$$
Therefore
$
ww^*\simeq \begin{pmatrix}1_k&0\\0&0 \end{pmatrix}
$, and hence
$$
yy^*\simeq \begin{pmatrix}0&0\\0&1_{N-k} \end{pmatrix},
\qquad
wtw^*\simeq \begin{pmatrix}t&0\\0&0 \end{pmatrix}
$$
We come to the desired expression.
\hfill $\square$


\section{Double cosets for infinite symmetric group}

\COUNTERS

{\bf\punct Infinite symmetric group.} Denote by $\S(\infty)$
the group of finitely supported
  permutations of a countable set $\{1,2,3,\dots\}$. We realize $\S(\infty)$
   as the group of infinite
0-1 matrices.
Next, we consider the group $\bfG:=\S(\alpha+m\infty)$,
which is the group $\S(\infty)$ represented  as a group of
block 0-1 matrices of the form (\ref{eq:g-big}).
Consider its subgroup $\bfK\simeq\S(\infty)$ realized
as the group of 0-1 matrices of the form (\ref{eq:u-big}).
We define the multiplication on
$\bfK\setminus \bfG/\bfK$ by the same formula
(\ref{eq:product-big}).

\sm

{\bf\punct The limit theorem.} We consider a group $G_N=\S(\alpha+m(k+N))$
represented as a subgroup of 0-1 matrices in $\U(\alpha+m(k+N))$.
Consider its subgroups $\S(\alpha+mk)$  and $K_N=\S(m(k+N))$, equip
$K_N$ with the uniform probability distribution $\kappa_N$.
For $g$, $h\in \S(\alpha+mk)$ we define the measure $\tau_{g,h}$
on $G_N$ as above (\ref{eq:tau-1})-(\ref{eq:tau-2}).
Also, we define the product $g\circ_N h$
$$
\S(\alpha+mk)\times \S(\alpha+mk)\to K_N\setminus G_N/K_N
.$$
as above (\ref{eq:gJh}).

\begin{theorem}
\label{th:4}
Fix $\alpha$. For any $k$ for any $\delta>0$ there exists
$N_0$ such that for all $N\ge N_0$ the measure of
the double coset $g\circ_N h$ is $>1-\delta$.
\end{theorem} 

{\sc Proof} is the same as for Theorems \ref{th:1}, \ref{th:2}, but it is more simple.
We can set $u=0$ outside a set of small measure. After this, all '$\simeq$' transform
to '='.
\hfill $\square$

\sm

{\sc Remark.} We observe that Theorem \ref{th:4} has a stronger form
than Theorems \ref{th:1} and \ref{th:2}. This form immediately  implies the multiplicativity
(as Theorem \ref{th:multi-1}.a)
theorem, see  Olshanski \cite{Olsh-topics}.

{\tt Math.Dept., University of Vienna,

 Nordbergstrasse, 15,
Vienna, Austria

\&

Institute for Theoretical and Experimental Physics,

Bolshaya Cheremushkinskaya, 25, Moscow 117259,
Russia

\&

Mech.Math. Dept., Moscow State University,
Vorob'evy Gory, Moscow

e-mail: neretin(at) mccme.ru

URL:www.mat.univie.ac.at/$\sim$neretin

wwwth.itep.ru/$\sim$neretin}

\end{document}